# PROPERTIES OF HIGHER CRITICISM UNDER STRONG DEPENDENCE

By Peter Hall and Jiashun Jin[1]

*University of Melbourne, University of California, Davis, and Purdue University*

The problem of signal detection using sparse, faint information is closely related to a variety of contemporary statistical problems, including the control of false-discovery rate, and classification using very high-dimensional data. Each problem can be solved by conducting a large number of simultaneous hypothesis tests, the properties of which are readily accessed under the assumption of independence. In this paper we address the case of dependent data, in the context of higher criticism methods for signal detection. Short-range dependence has no first-order impact on performance, but the situation changes dramatically under strong dependence. There, although higher criticism can continue to perform well, it can be bettered using methods based on differences of signal values or on the maximum of the data. The relatively inferior performance of higher criticism in such cases can be explained in terms of the fact that, under strong dependence, the higher criticism statistic behaves as though the data were partitioned into very large blocks, with all but a single representative of each block being eliminated from the dataset.

## 1. Introduction.

1.1. *Decision-making using sparse, faint information.* Modern data acquisition routinely produces massive, complex datasets in many scientific areas, for example, genomics, astronomy and functional MRI. The need for fast and effective analysis in these settings poses challenging statistical problems, one of which is how to reliably detect the presence of a sparse, faint signal. Here, data are available on a large number of observation units (or hypothesis tests, or transform coefficients, etc.), which may or may not contain

Received December 2006; revised April 2007.
[1]Supported in part by NSF Grants DMS-05-05423 and DMS-06-39980.

a signal; the signal, when present, is faint and is dispersed across different observation units (e.g., vector components) in an unknown fashion.

The situation arises naturally in a range of application areas, for example, early detection of airborne bio-terror and syndromic surveillance (Donoho and Jin [7] and Anon [1]), covert communications (Donoho and Jin [7]) and non-Gaussian detection of the cosmic microwave background (CMB) (Cayon, Jin and Treaster [5] and Jin et al. [20]). In these examples the desired signal may represent the outbreak of certain disease, or covertly attached signals in a white noise channel, or non-Gaussian signatures in the CMB, and so on. In such cases the signal is either highly sparse, or faint in its individual presences, or both. The sparsity and faintness play an intricate duet which calls for nontraditional methods for signal detection.

The signal detection problem is closely connected to that of multiple hypothesis testing. Indeed, if in the former setting we associate each observation with a null hypothesis, which asserts that the signal is not present and the observation is pure noise, then the signal detection problem can be viewed as attempting to determine whether at least one of the null hypotheses is false. From this point of view, related work has been done in, for example, problems of assessing the accuracy of random number generators; see, for example, Knuth [22]. Review-type accounts of multiple hypothesis testing include those of Hochberg and Tamhane [13], Pigeot [25], Dudoit, Shaffer and Boldrich [9], Bernhard, Klein and Hommel [3] and Lehmann and Romano ([23], Chapter 9).

Of course, the problem of (1) discovering which null hypotheses are false is more difficult than that of (2) estimating the proportion of false null hypotheses, which in turn is more challenging than (3) determining whether at least one null hypothesis is false. Work on these respective problems includes that of (1) Benjamini and Hochberg [2], Genovese and Wasserman [11], Donoho and Jin [8], Efron et al. [10] and Storey, Dai and Leek [26]; (2) Swanepoel [27], Efron et al. [10], Cai, Jin and Low [4], Genovese and Wasserman [11], Storey, Dai and Leek [26], Jin [18], Jin and Cai [19], Jin, Peng and Wang [21] and Meinshausen and Rice [24]; (3) Donoho and Jin [7], Delaigle and Hall [6], Hall, Pittelkov and Ghosh [12], Ingster [14, 15], Jin [17] and Jager and Wellner [16].

1.2. *Higher criticism methods for independent data.* Inspired by ideas of Tukey [28], Donoho and Jin [7] proposed higher criticism methods for signal detection in the presence of a white noise. The technique is based on assessing the statistical significance of the number of significant results in a long sequence of hypothesis tests, which can be either formal or informal. The principle has found a variety of applications, for example, to non-Gaussian detection (Cayon, Jin and Treaster [5] and Jin et al. [20]), goodness of fit (Jager and Wellner [16]) and classification (Hall, Pittelkov and Ghosh [12] and Delaigle and Hall [6]).

A higher criticism statistic, computed from standard normal data $X_i$ to which might be added sparsely distributed, positive signals, can be defined by

$$\mathrm{hc}_n = \sup_{t \in \mathcal{S}_n} \frac{\sum_{1 \le i \le n} \{I(X_i > t) - \bar{\Phi}(t)\}}{\{n\Phi(t)\bar{\Phi}(t)\}^{1/2}}, \tag{1.1}$$

where $\mathcal{S}_n$ denotes a subset of the real line, $\Phi$ is the standard normal distribution function and $\bar{\Phi} = 1 - \Phi$. If positive signals are added to the $X_i$'s, then the numerator on the right-hand side of (1.1) takes relatively large values. Provided the size and number of added signals are sufficiently great, $\mathrm{hc}_n$ exceeds, with high probability, a critical point calculated under the assumption of "no signal." The size of this exceedance can be used as a basis for detecting the presence of the signals.

1.3. *The case of dependent data.* It can be shown (see Delaigle and Hall [6]) that the main features of higher criticism do not alter under conditions of short-range dependence, for example, if the data $X_i$ come from a moving-average process with exponentially decaying weights. However, the nature of the higher criticism statistic can change dramatically under strong dependence. Provided an appropriate critical point is employed to assess significance, higher criticism can still be used effectively in such cases, although its performance does not necessarily compare well with that of competing, difference-based signal detectors. Nevertheless, in order to be effective the latter approaches can require significant information to be known about the signal, and so their theoretical attractiveness may not necessarily be evidenced in practice.

To explore and elucidate these issues we shall treat dependent data generated under a simple autocovariance model:

$$\rho_n(k) = \max(0, 1 - |k|^\alpha \ell_n^{-\alpha}), \tag{1.2}$$

where $\alpha > 0$ and $\ell_n$ denotes a positive sequence diverging to infinity. If $X_i$ is a zero-mean, stationary Gaussian process with $\mathrm{cov}(X_{i_1}, X_{i_2}) = \rho_n(i_1 - i_2)$, then data values lagged $\ell_n$ or more apart are independent. Therefore, by choosing $\ell_n$ to diverge more slowly we reduce the range of dependence. Motivation for contexts such as this, and for the dependent data case more generally, is given in Section 1.4.

We shall show that for large $n$, and to a first approximation, the numerator at (1.1) behaves like $\ell_n \sum_j \{I(Z_j > t) - \bar{\Phi}(t)\}$, where there are just $n/\ell_n$ terms in this series and the $Z_j$'s are independent and normal $\mathrm{N}(0,1)$. That is, the higher criticism statistic behaves as though the data were partitioned into $\ell_n$ blocks, and were identically equal to $Z_j$ within the $j$th block. (Here it is assumed that $1 \le \ell_n \le n$.) Reflecting this property, the number of signals present should be increased by the factor $\ell_n$ if performance is to be similar to that in the case of independence; we shall make this precise in Section 3.4. This blockwise view provides considerable insight into properties of the higher criticism statistic, although it conceals the fact that, under strong dependence, the differences $X_{i+1} - X_i$ can be used to construct powerful tests.

1.4. *Speckle imaging.* One of the challenging aspects of earth-based astronomical imaging is removing the effects of atmospheric turbulence. The atmosphere is in constant motion, and it alters a point-like sharp signal to a blurred signal whose position changes, sometimes rapidly, over time. Speckle imaging involves creating many images, one after the other (e.g., thousands in an hour are possible using current technology), and combining these rapid "snapshots" to gain greater information (e.g., higher resolution) than could be managed using conventional imaging. A disadvantage of this approach is that it gives images with relatively low brightness, and so many short-exposure images are "stacked" to achieve greater sensitivity. There is a variety of ways of aligning stacked images, for example, by using the brightest point (a "speckle") as the locus.

Of course, the brightest speckle moves relative to the point source in the heavens, for example, because of atmospheric effects and telescope vibration. Moreover, it can subdivide into more than one speckle, due to atmospherical effects. Therefore image stacking retains significant levels of noise, although it produces images which enjoy relatively high signal-to-noise ratios at high frequencies. While image stacking involves spatial, two-dimensional data, intellectually the issues are most transparently tackled in a one-dimensional setting, without making intrinsic changes. In this paper we analyze a one-dimensional model for the type of data obtained from image stacking, and report on the implications of the model.

The stacking procedure results, first, in very faint point sources, often only a pixel or two wide and, in many cases, representing the same heavenly point source but located in different (perhaps multiple) positions in each image; and, second, in a background noise process that is defined essentially in the continuum, and whose correlation can extend over many pixels, especially if pixel width is small. In particular, the correlation between the noise at pixels $i$ and $j$ can be represented as $1 - g(|i-j|/n)$, where $n^{-1}$ denotes pixel width and $g$ is a smooth function. The autocovariance model at (1.2) is an idealized form of this.

As technology improves, the distance between adjacent pixels becomes smaller (equivalently, the value of $n$ becomes larger), and the number of images in a stack increases (implying that the interpolated continuum model in the previous paragraph applies more closely). Work of Jin [17] uses higher criticism methods in this setting, but against the background (and the independence assumptions) of Donoho and Jin [7]. The present paper takes a more careful look at what can occur in the case of dependence, and in particular reveals the suboptimality of higher criticism there.

To extend these concepts to other settings, suppose we have data from a process $Z$ which is a discretization (e.g., on a pixel grid) of a continuum process for which the autocorrelation function is $1-g$, where $g$ is smooth. As the discretization becomes finer, the strength of dependence of the discrete process increases, and the results developed in this paper become relevant. In the astronomical example, greater fineness is the result of improved imaging methods, but in other settings the cause may be different.

In practice, a reasonable amount of information is known about the signal, and as in all imaging problems, further advice on choice of the critical point can be gained by visual experimentation. (Sometimes the differences between imaging problems and curve estimation are overlooked in this respect.)

## 2. Overview.

2.1. *Time-series model, and its immediate consequences.* Let $X_{n1}, X_{n2}, \ldots$ denote a stationary Gaussian process, with zero mean and autocovariance given by a slight specialization of (1.2):

$$\rho_n(k) = \max(0, 1 - |k|^{\alpha} n^{-\alpha_0}), \tag{2.1}$$

where $\alpha, \alpha_0 > 0$. That is, $\operatorname{cov}(X_{ni_1}, X_{ni_2}) = \rho_n(i_1 - i_2)$ for each pair $(i_1, i_2)$. The length of the range of dependence of this process increases with both $\alpha_0$ and $\alpha^{-1}$.

Let $\kappa = \alpha_0/\alpha$. Then $\ell_n$, in (1.2), is effectively equal to $n^{\kappa}$. If we are considering crossings of levels proportional to $\sqrt{\log n}$, then, to a first approximation, there is a high probability that any crossing is succeeded by "almost" $n^{\kappa}$ further such crossings, although not by $n^{\kappa}$ crossings. To make this claim more concise we note that if $t = \sqrt{2q \log n}$, where $q > 0$, then

$$\text{for each } \varepsilon > 0, \qquad P(X_{ni} > t \text{ for } 2 \le i \le n^{\kappa - \varepsilon} | X_{n1} > t) \to 1. \tag{2.2}$$

On the other hand,

$$P(X_{ni} > t \text{ for } 2 \le i \le n^{\kappa} | X_{n1} > t) \to 0. \tag{2.3}$$

Result (2.2), under the side condition $\alpha \ge 2$, follows from Theorem 3.2 in Section 3. Property (2.3) is simpler to derive; note that if $X$ and $Y$ are jointly normally distributed with zero means, unit variances and covariance $\rho$, then $P(X > t | Y > t) \to 0$ as $t \to \infty$.

2.2. *Blockwise decomposition of higher criticism statistic.* Suppose we observe only the first $n$ terms, that is, $X_{n1}, \ldots, X_{nn}$, in the time series. Then (2.2) and (2.3) indicate that, if we are addressing level-crossings on a $\sqrt{\log n}$ scale, the sum of indicator functions that defines the higher criticism statistic can be written approximately as a sum of $\max(1, n^{1-\kappa})$ indicators of time-series variables lagged $\min(n, n^{\kappa})$ apart:

$$\sum_{i=1}^{n} I(X_{ni} > t) \approx \min(n, n^{\kappa}) \sum_{j=0}^{\max(1,n^{1-\kappa})} I(X_{n,jn^{\kappa}+1} > t). \tag{2.4}$$

Here we interpret $jn^{\kappa}$ as $j\lfloor n^{\kappa} \rfloor$, where $\lfloor x \rfloor$ denotes the least integer not strictly less than $x$. We may consider $X_{n,jn^{\kappa}+1}$, appearing in the argument of the indicator variable on the right-hand side of (2.4), as representing the $j$th "block" of time-series values, that is, as representing the data $X_{n,jn^{\kappa}+1}, X_{n,jn^{\kappa}+2}, \ldots, X_{n,(j+1)n^{\kappa}}$.

The indicator variables on the right-hand side of (2.4) are independent and identically distributed as normal $\mathrm{N}(0, 1)$, and so (2.4) implies that the higher criticism statistic is approximately equal to its version for the smaller, "effective sample size" of $\max(1, n^{1-\kappa})$, subsequently multiplied by the factor $\min(n, n^{\kappa})$. This property gives insight, discussed in Section 2.3, into higher criticism for strongly dependent data.

Even if block size is reduced from $n^{\kappa}$ to $n^{\kappa-\varepsilon}$, for some $\varepsilon > 0$, then it is unreasonable to expect blockwise decompositions to hold uniformly in $t$. To appreciate why, consider slowly increasing the value of $t$ until a point $t'$, say, is reached where, for a particular index $j$, at least one of the variables in the $j$th block first fails to exceed the level $t$. Although the variables in the $j$th block are highly correlated, they have a proper joint distribution. Therefore, as we increase $t$ beyond $t'$, the other variables in the block will fail one by one to exceed $t$. They will not fail simultaneously, although, since the correlation is high, they will, with high probability, all fail within a relatively small interval of values of $t$.

2.3. *Summary of properties of higher criticism statistic.* We shall argue that behavior of the higher criticism statistic can be decomposed into two cases, "degenerate" and "nondegenerate." The degenerate case arises when $\kappa \geq 1$. Here, in view of (2.2), the effective sample size, for crossings of a level on a $\sqrt{\log n}$ scale, is just 1. Therefore, if $\kappa \geq 1$, then the strongly dependent nature of the data effectively restricts us, when using the conventional higher criticism statistic, to working with a single data value.

If $\kappa < 1$, then the problem is nondegenerate, and (2.2) and (2.3) imply that the effective sample size is $N = n^{1-\kappa}$. In this case, if $v > 1$, then the probability that an exceedance of the level $\sqrt{2v \log N}$ occurs, among the $N$ independent data, converges to zero, and so we should confine attention to levels for which $v < 1$. For simplicity, we skip discussion of the case of $v = 1$. Taking $\sqrt{2q \log n}$ to be the level of the signal, and equating $\sqrt{2v \log N}$ to $\sqrt{2q \log n}$, we see that "$v < 1$" translates to the condition "$q < 1 - \kappa$." Therefore, we argue:

(2.5) CLAIM. In keeping with the approach developed for independent data, when $\kappa < 1$ the higher criticism statistic should be used to address crossings of levels no higher than $\sqrt{2q \log n}$, where $q < 1 - \kappa$.

Theorem 3.4 will justify the claim.

In the discussion above we treated the "null" setting, where no signal is present. In conventional higher criticism (Donoho and Jin [7]), the nonnull case is constructed by distributing $n^{1-\beta}$ signals, each equal to $\sqrt{2r \log n}$ for some $r \in (0, 1)$, independently and uniformly among the $n$ noise variables $X_i$. [To make the signal-detection problem nontrivial we take $\beta \in (\frac{1}{2}, 1)$.] To keep faith with the blockwise treatment suggested in Section 2.2, we should ideally distribute $N^{1-\beta}$ blocks of signals, each block comprising $n^\kappa$ signals and each signal equal to $\sqrt{2r \log N}$ where $0 < r < 1$, among the $N$ blocks. Of course, this does not happen; the blocks are fictions of our mathematical argument, and are not respected by the physical process that produces the signal. Nevertheless, as we shall show in Section 3.4, a close parallel with the case of independent data emerges if we add $N^{1-\beta} n^\kappa$ signals, distributed among the $n$ time-series values $X_{ni}$.

It follows that, if $\kappa < 1$ and we distribute $N^{1-\beta} n^\kappa$ signals randomly and uniformly among the $n$ time-series values, then the higher criticism statistic for the time-series dataset $X_{n1}, \dots, X_{nn}$ is well approximated by a constant multiple of its counterpart when the time-series is sampled only $N$ times, once at every $n^\kappa$ points, and the $N^{1-\beta}$ signals are distributed among the $N$ sampled points. The resulting subseries is comprised only of independent data.

In the case of a strongly dependent time-series $X_{ni}$, even if the magnitude of the signal is as small as $n^{-C_1}$ for $C_1 > 0$ not too large, the presence of the signal can be determined accurately merely by deciding that the signal is present if, for some $i \in [1, n-1]$, $|X_{n,i+1} - X_{ni}| > n^{-C_2}$, where $0 < C_2 < C_1$. For this simple signal detector, the probability that the signal is not detected when it is not present, and the probability that it is detected when it is present, both converge to 1. However, effectiveness of the method requires information about the range of dependence. Details will be given in Theorem 3.8.

2.4. *Signal detection using the maximum.* A commonly used statistic for signal detection is the maximum of the observed data, $\text{Max}_n = \max(X_{n1}, \dots, X_{nn})$. In the case of independent, $N(0, 1)$ noise, the signal is deemed to be present if $\text{Max}_n > \sqrt{2 \log n}$, and not present otherwise; see, for example,

Donoho and Jin [7]. This "max classifier" has less sensitivity than higher criticism in the white noise setting, but is more robust to dependence. Details will be given in Section 3.4.

## 3. Main results.

3.1. *Variance of numerator of higher criticism statistic.* Recall that $\kappa = \alpha_0/\alpha$, where $\alpha, \alpha_0 > 0$ are parameters governing the autocovariance model at (2.1). Given two functions $a_1$ and $a_2$ satisfying $0 < a_1(t) \le a_2(t) < \infty$, write $\langle a_1(t), a_2(t)\rangle$ to denote a quantity which, uniformly in $n \ge 1$ and $t \ge 1$, is bounded between $C_1 a_1(t)$ and $C_2 a_2(t)$, for constants $C_1, C_2 > 0$. Let $t_1 = |t| + 1$, for arbitrary real $t$.

Our first theorem describes the variance of the argument of the higher criticism statistic, uniformly on the positive half-line.

THEOREM 3.1. *If the autocovariance of the time-series $X_{ni}$ is given by (2.1), then there exist constants $B_1, B_2 \ge 1$ such that, uniformly in $n \ge 1$ and all $t$,*

$$\operatorname{var}\left\{\sum_{i=1}^{n} I(X_{ni} > t)\right\} = \langle t_1^{-B_1}, t_1^{B_2}\rangle n^{\min(\kappa+1,2)} e^{-t^2/2}. \tag{3.1}$$

REMARK 3.1 (*Impact of blockwise properties on variance*). Recall, from (2.4), that the numerator of the higher criticism statistic can be approximated by a construction where the terms in the numerator are grouped into blocks of length $n^\kappa$, and the indicator functions that represent respective blocks are independent. Referring to (2.4), and writing simply $X$ for a random variable with the standard normal distribution, the variance of this approximation can be seen to equal

$$\begin{aligned} &\max(1, n^{1-\kappa}) \operatorname{var}\{\min(n, n^\kappa) I(X > t)\} \\ &\quad = n^{\min(\kappa+1,2)} \Phi(t)\{1 - \Phi(t)\} = \langle t_1^{-1}, t_1^{-1}\rangle n^{\min(\kappa+1,2)} e^{-t^2/2}. \end{aligned} \tag{3.2}$$

The right-hand sides of the specific formula (3.1), and its approximation (3.2), are close to one another.

3.2. *Blockwise properties of higher criticism statistic.* Here we state results that underpin (2.2), (2.4) and (2.5). Recall that $\kappa = \alpha_0/\alpha$.

THEOREM 3.2. *Assume that the stationary, zero-mean Gaussian process $X_{ni}$ has autocovariance given by (2.1), and that $0 < \lambda < \kappa$. Then, for each $\eta > 0$,*

$$P\{I(X_{ni} > t) = I(X_{n1} > t) \text{ for } 1 \le i \le n^\lambda\} = 1 - O(n^{\lambda - (\alpha_0/2) + \eta} e^{-t^2/2}), \tag{3.3}$$

*uniformly in all $t$.*

Since $\lambda < \kappa$ can be chosen arbitrarily close to $\kappa$, then Theorem 3.2, and related properties such as those at (2.3) and in Theorem 3.3, demonstrate that the process of indicator values, $I(X_{ni} > t)$, can be divided approximately into blocks of length $n^\kappa$. Here the approximation is logarithmic: the logarithm of block length is approximately equal to $\log(n^\kappa)$. The difference in size of the higher criticism statistic, in the presence and in the absence of signals, respectively, is small on the same scale; the change in size is by a factor of $n^\varepsilon$, where $\varepsilon$ can be very small although always strictly positive, and depends on the distance that the signal lies above the detection boundary. See Theorem 3.7. Therefore, measurement of block size on a logarithmic scale is appropriate in the present setting.

A corollary of Theorem 3.2 is that, if $\alpha \ge 2$, which ensures that $\kappa - (\alpha_0/2) \le 0$; and if $t_n$ is any sequence of positive constants such that $t_n/n^\eta \to 0$ for each $\eta > 0$; then:

$$\text{for each } \varepsilon > 0, \qquad \sup_{t:|t|\le t_n} |1 - P(X_{ni} > t \text{ for } 2 \le i \le n^{\kappa-\varepsilon} | X_{n1} > t)| \to 0. \tag{3.4}$$

This is a strong form of (2.2).

Result (3.4) also implies that if $\kappa > 1$, then, with high probability, either all the data in the sample $X_{n1}, \ldots, X_{nn}$ are above the level $t$, or all are below that level. Indeed, with $t_n$ as before,

$$\sup_{t:|t|\le t_n} |1 - P(X_{ni} > t \text{ for } 2 \le i \le n | X_{n1} > t)| \to 0. \tag{3.5}$$

However, as (2.3) indicates, (3.5) fails when $\kappa = 1$.

Our next result quantifies the approximation at (2.4). We already know, from (3.5), that if $\kappa > 1$, then (2.4) holds in the sense that, for $J = 0$ or 1,

$$P\left\{\sum_{i=1}^n I(X_{ni} > t) = nI(X_{n1} > t) | I(X_{n1} > t) = J\right\} \to 1, \tag{3.6}$$

uniformly in $|t| \le t_n$. For $\kappa \le 1$ the approximation is a little more tricky, in that, as indicated at the end of the previous paragraph, block length cannot be quite as long as $n^\kappa$ if the sum of indicator functions, on the left-hand side of (2.4), is to decompose into blocks with high probability. Moreover, when there is more than one block it is awkward to sharpen the result by conditioning, as at (3.6). Theorem 3.3, below, gives a concise account of subdivision into blocks of length $\lfloor n^\lambda \rfloor$, where $\lambda < \kappa$ can be chosen arbitrarily close to $\kappa$.

Before stating Theorem 3.2 we give a little notation. Recall that $\kappa = \alpha_0/\alpha$. If $\kappa > 1$, take $b = 1$ and $\mathcal{B}_1 = \{1, \ldots, n\}$. If $\kappa \le 1$ and $\lambda \in (0, \kappa)$, partition the integers $1, \ldots, n$ into $b$ consecutive, adjacent blocks $\mathcal{B}_1, \ldots, \mathcal{B}_b$, where the first $b-1$ blocks contain just $\lfloor n^\lambda \rfloor$ integers and the last block is of length between 1 and $\lfloor n^\lambda \rfloor$. (Thus, $b \sim n^{1-\lambda}$.) Let $s_j$ denote the least element of $\mathcal{B}_j$.

THEOREM 3.3. *If the assumptions in Theorem 3.2 hold, then, for each* $\eta > 0$,

$$P\{I(X_{ni} > t) = I(X_{ns_j} > t) \text{ for } i \in \mathcal{B}_j \text{ and } 1 \le j \le b\} \\ = 1 - O(n^{1-(\alpha_0/2)+\eta} e^{-t^2/2}) \tag{3.7}$$

*uniformly in all* $t$.

REMARK 3.2 (*Blockwise clumping of indicator variables*). If we add to the assumptions of Theorem 3.3 the condition $\alpha_0 > 0$, then (3.7) implies that

$$\sup_{-\infty < t < \infty} |1 - P\{I(X_{ni} > t) = I(X_{ns_j} > t) \text{ for } i \in \mathcal{B}_j \text{ and } 1 \le j \le b\}| \to 0$$

uniformly in all $t$. This result, and (3.4), underpin the blockwise decomposition of the higher criticism statistic, discussed in Section 2.2.

Finally in this section we justify (2.5). We argue that if $t = \sqrt{2q \log n}$, where $q > 1 - \kappa$, then, in the majority of samples, none of the noise data exceed the level $t$. See (3.8) below. Therefore, if the signal is at level $t$ or larger, and if it is added to the noise, then it will be very easy to detect. Hence, to make the signal-detection problem reasonably difficult, we should ensure that $q < 1 - \kappa$, as claimed at (2.5).

THEOREM 3.4. *Assume that the stationary, zero-mean Gaussian process* $X_{ni}$ *has autocovariance given by (2.1), and that* $\kappa < 1$. *Then, if* $t = \sqrt{2q \log n}$ *where* $q > 1 - \kappa$,

$$P(X_{ni} > t \text{ for some } 1 \le i \le n) \to 0. \tag{3.8}$$

3.3. *Size of higher criticism statistic.* Here we address properties of $\mathrm{hc}_n$, defined at (1.1), in cases where $\kappa < 1$. Our analysis of blockwise characteristics of the higher criticism statistic has already shown that the case $\kappa > 1$ is relatively uninteresting, since there, almost all the data $X_{n1}, \ldots, X_{nn}$ cross a given level $t$ at the same time.

We first address the nature of the set $\mathcal{S}_n$ in the definition of $\mathrm{hc}_n$. If the data $X_{ni}$ were independent and identically distributed, or equivalently, if $\kappa = 0$, then we could take $\mathcal{S}_n = [-t_n, t_n]$ where $t_n = \sqrt{2s \log n}$ and $0 < s < 1$. Taking $s \ge 1$ is inappropriate, since the version of Theorem 3.4 that applies in the setting of independent data then implies that the level $t_n$ is hardly ever exceeded by the data in a sample of $n$ standard normal random variables; the level is too large.

To remove this difficulty we should instead consider a value of $t_n$ such that $n\bar{\Phi}(t_n)$ is bounded away from zero. A borderline sequence of this type is $t_n = \{2\log n - \log(C\log n)\}^{1/2}$, where $C$ is any positive constant. Here, $n\bar{\Phi}(t_n)$ converges to $(C/\pi)^{1/2}$. We know from the results discussed in Section 3.2 that, when $\kappa < 1$, the data behave as though they were partitioned into $n^{1-\kappa}$ blocks, where the indicators $I(X_{ni} > t)$ are identical within the block. Therefore, when $\kappa < 1$ we should ask that $n^{1-\kappa}\bar{\Phi}(t_n)$, rather than $n\bar{\Phi}(t_n)$, be bounded away from zero. This appreciation motivates the assumption on $t_n$ imposed in the theorem below.

THEOREM 3.5. *Assume that the stationary, zero-mean Gaussian process $X_{ni}$ has autocovariance given by (2.1), that $\kappa < 1$, and that $t_n \to \infty$ in such a manner that $n^{1-\kappa}\bar{\Phi}(t_n)$ is bounded away from zero. Then, for all $\eta > 0$,*

$$P\left[\sup_{|t|\le t_n}\left|\frac{\sum_i\{I(X_{ni}>t)-P(X_{ni}>t)\}}{\{n\Phi(t)\bar{\Phi}(t)\}^{1/2}}\right| > n^{(\kappa+\eta)/2}\right] \to 0. \tag{3.9}$$

Theorem 3.5 points to the size of critical point appropriate for a test of significance involving the higher criticism statistic, as follows. A multivariate central limit theorem for values of $\sum_i I(X_{ni} > t)$, for a fixed but arbitrarily large number of different $t$'s, implies that if the critical point $c_n(\alpha)$ is to satisfy

$$P\left[\sup_{|t|\le t_n}\frac{\sum_i\{I(X_{ni}>t)-P(X_{ni}>t)\}}{\{n\Phi(t)\bar{\Phi}(t)\}^{1/2}} > c_n(\alpha)\right] = \alpha, \tag{3.10}$$

where $0 < \alpha < 1$ is fixed, then $c_n(\alpha) = n^{\kappa/2}d_n(\alpha)$, where $d_n(\alpha)$ diverges to infinity. (The central limit theorem can be proved using the method of moments.) On the other hand, Theorem 3.5 shows that $d_n(\alpha)$ is no larger than $O(n^\eta)$, for any $\eta > 0$. These considerations lead to the following corollary to Theorem 3.5.

THEOREM 3.6. *If $X_{ni}$ is a stationary, zero-mean Gaussian process with autocovariance given by (2.1), if $0 < \kappa < 1$, if $t_n \to \infty$ such that $n^{1-\kappa}\bar{\Phi}(t_n)$ is bounded away from zero, and if $c_n(\alpha)$ is defined by (3.10), then $c_n(\alpha) = n^{\kappa/2}d_n(\alpha)$ where, as $n$ increases, $d_n(\alpha)$ diverges to infinity for each fixed $\alpha$, but equals $O(n^\eta)$ for each $\eta > 0$.*

REMARK 3.3 (*Calibration*). In practice, a critical point would be determined either by experience with the time-series $X_{ni}$, or by simulation from a model for those data. In either case, $c_n(\alpha)$ would, in effect, be found empirically. However, as Theorem 3.7 will show, it is not important that $c_n(\alpha)$ be determined particularly accurately.

Let $0 < \kappa < 1$ and $N = n^{1-\kappa}$, as in Section 2, and consider the impact of adding a signal, equal to

$$\nu = \sqrt{2r \log N}, \tag{3.11}$$

to just $N^{1-\beta} n^{\kappa}$ of the standard normally distributed data $X_{n1}, \ldots, X_{nn}$. This changes the data from $X_{ni}$ to

$$Y_{ni} = X_{ni} + I_{ni}\nu, \tag{3.12}$$

where $I_{ni}$ is a process of zeros and ones. This suggests that we add a signal to all $n^{\kappa}$ time-series values in each of $N^{1-\beta}$ blocks. However, we do not need to add the signals in this blockwise way. They can be added in any deterministic manner, or in any random way that is stochastically independent of the time-series $X_{ni}$. We take $\frac{1}{2} < \beta < 1$, to make the signal-detection problem relatively difficult.

The detection boundary of Donoho and Jin [7] (see also Ingster [14, 15] and Jin [17]) is the locus of points $(\beta, r)$, with $\frac{1}{2} < \beta < 1$, such that

$$r = \begin{cases} \beta - \frac{1}{2}, & \text{if } \frac{1}{2} < \beta \le \frac{3}{4}, \\ (1 - \sqrt{1-\beta})^2, & \text{if } \frac{3}{4} < \beta < 1. \end{cases} \tag{3.13}$$

The theorem below shows that, if $c_n(\alpha)$ is the critical point defined by (3.10); and if we assert that the signal is present if the higher criticism statistic exceeds $c_n(\alpha)$, or even if it exceeds a bound that is larger than that quantity to a small but fixed polynomial extent; then the probability that we make the correct decision when the signal is present, converges to 1 if $(\beta, r)$ lies above the boundary. Below we take $N = n^{1-\kappa}$, where $0 < \kappa < 1$, and write $\mathrm{hc}_n^{\mathrm{sig}}$ for the version of $\mathrm{hc}_n$ when $X_{ni}$, at (1.1), is replaced by $Y_{ni}$.

THEOREM 3.7. *Let $X_{ni}$ be a stationary Gaussian process with zero mean and autocovariance as at (2.1), and let $Y_{ni}$ be as at (3.12), with the $I_{ni}$'s independent of $X_{n1}, \ldots, X_{nn}$ and just $N^{1-\beta} n^{\kappa}$ of them equal to 1. Assume that $n^{1-\kappa}\bar{\Phi}(t_n)$ is bounded away from zero, and $(\beta, r)$ lies strictly above the detection boundary. Then, for $\delta > 0$ sufficiently small,*

$$P(\mathrm{hc}_n^{\mathrm{sig}} > n^{(\kappa/2)+\delta}) \to 1. \tag{3.14}$$

Define $\xi = (1-\beta)(1-\kappa)$. If we consider that $n^{1-\beta'} = N^{1-\beta} n^{\kappa} = n^{\xi+\kappa}$ signals, each of size $\sqrt{2r' \log n} = \sqrt{2r \log N}$, have been added to the time-series $X_{n1}, \ldots, X_{nn}$, then $(\beta', r')$ is related to $(\beta, r)$ by the formulas $\beta' = \beta(1-\kappa)$ and $r' = r(1-\kappa)$. Figure 1 graphs the detection boundary determined by (3.13), in terms of $(\beta', r')$ rather than $(\beta, r)$.

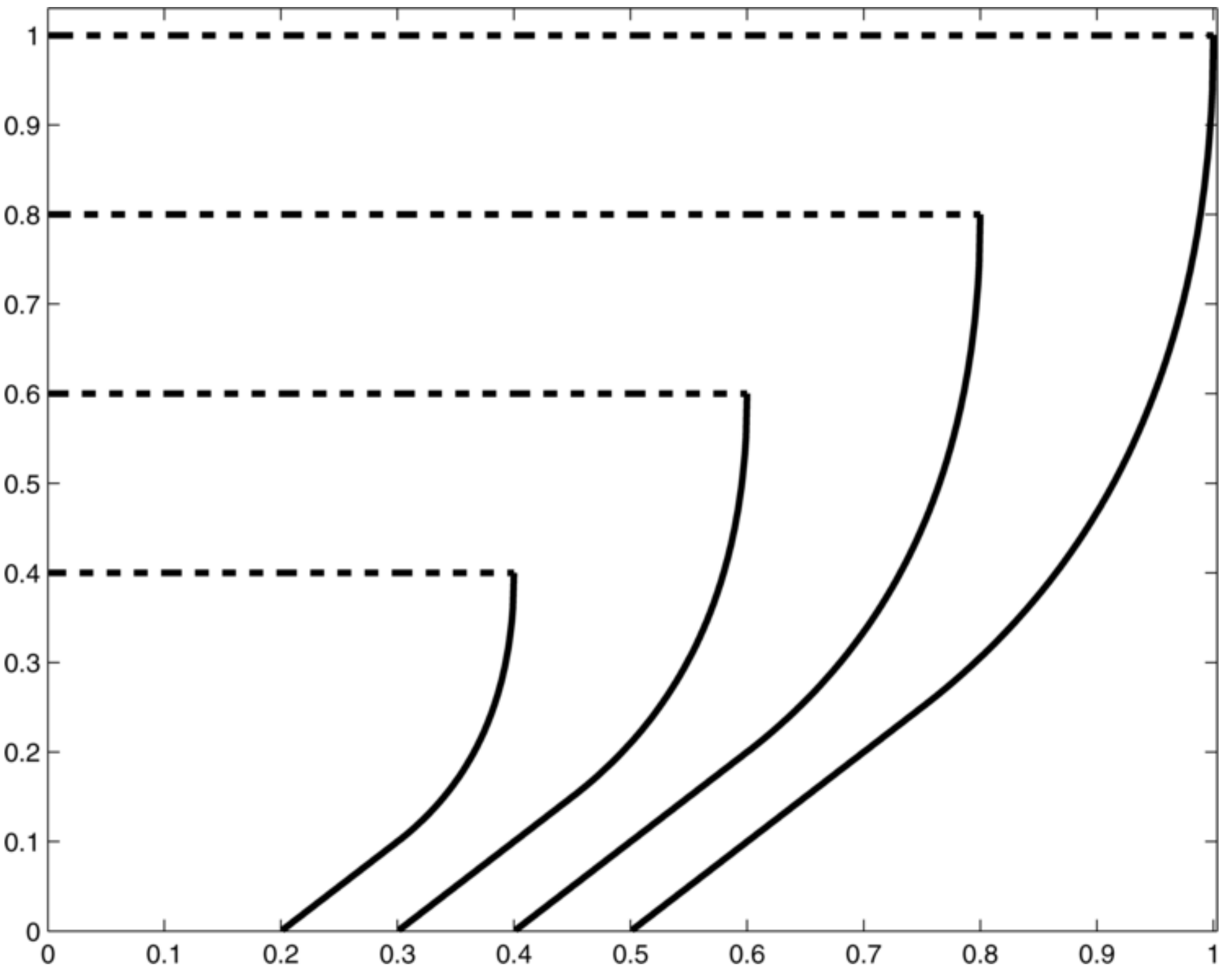

Fig. 1. *Detection boundary graphed in terms of $(\beta', r')$ for different $\kappa$. Recall that $\beta' = (1-\kappa)\beta$ and $r' = (1-\kappa)r$. Each curve depicts the boundary shown at (3.14), with variables $(\beta, r)$ rescaled to $(\beta', r')$. The curves, in order from top to bottom, are for the cases $\kappa = 0.6, 0.4, 0.2, 0$, the last denoting the context of independent data. The horizontal dashed lines show graphs of $r' = 1-\kappa$ for the four respective values of $\kappa$.*

Remark 3.4 [*Performance of higher criticism when $(\beta, r)$ lies below the detection boundary*]. An argument similar to that given in the proof of Theorem 3.7 may be used to prove that, if $(\beta, r)$ lies strictly below the detection boundary, and if the positive constants $d_n$ diverge at rate $n^\delta$ for each $\delta > 0$, although nevertheless sufficiently fast, then $P(\mathrm{hc}_n^{\mathrm{sig}} > n^{\kappa/2} d_n) \to 0$. Therefore, the higher criticism statistic cannot be relied on to give accurate detection when $(\beta, r)$ lies below the detection boundary.

If we are aware of the presence of strong dependence, and can exploit it through an accurate mathematical model for both the error process and the signal, then the dependence can be utilized to produce a signal detector with a high degree of sensitivity. To see why, consider adding a much smaller signal than before, this time of size only $\nu' = \pm 2(rn^{-\alpha_0} \log n)^{1/2}$ where $r > 1$, to between one and $n-1$ points of the time-series $X_{n1}, \ldots, X_{nn}$, thereby obtaining the data $Y'_{ni}$, say. Let $1 \le C < r$, and suppose that we determine the signal to be present if

$$\max_{1 \le i \le n} |Y'_{n,i+1} - Y'_{ni}| > 2(C n^{-\alpha_0} \log n)^{1/2}.$$

If this inequality fails, we determine the signal to be absent. We shall call this the "neighbor-difference detector" (NDD). In both conception and effect it is rather like assessing the time-series $Y'_{ni}$ using a high-pass filter.

Our next result shows that NDD enjoys a high degree of sensitivity. Depending on how much information is available about the signal and the noise process, NDD might require, relative to other methods, greater knowledge of the covariance at (2.1). The main issue in practice, however, is determining when the level of correlation has increased to such an extent that methods based on the assumption of independence are no longer competitive. As noted in Section 1.4, this point can be reached relatively subtly, for example, as imaging technology improves through decreases in pixel size.

THEOREM 3.8. *If $1 \le C < r$, then the probability that NDD correctly determines that a signal is present, given that it is, and the probability that the detector correctly determines that a signal is not present, given that it is not, both converge to 1 as $n \to \infty$.*

Note that, although Theorem 3.8 treats only one signal, that signal is permitted to be added at any number of points, between one and $n-1$, and can be polynomially small in size, rather than logarithmically large as would be required in the case of higher criticism. Therefore, when using NDD the signal can be much smaller in size, and much less frequently present, than in the case of higher criticism; and nevertheless NDD manages to detect it.

3.4. *Properties of* $\mathrm{Max}_n$. The max classifier was defined in Section 2.4. In the case of independent $\mathrm{N}(0,1)$ data, the rule given there leads to asymptotically correct classification if $r > (1 - \sqrt{1-\beta})^2$, and to correct classification with limiting probability $\frac{1}{2}$ if $r < (1 - \sqrt{1-\beta})^2$. The resulting detection boundary, with equation $r = (1 - \sqrt{1-\beta})^2$, coincides with that for higher criticism if $\frac{3}{4} \le \beta < 1$, but is above the higher criticism boundary if $\frac{1}{2} < \beta < \frac{3}{4}$. However, unlike higher criticism, the maximum-based detector maintains its performance in the presence of a polynomially large amount of dependence.

To appreciate this point, given $\beta \in (\frac{1}{2}, 1)$ and a nonnegative sequence $a_n$ decreasing to zero, let $\mathcal{X}(\beta, a_n)$ denote the set of all distributions of time-series $X_{n1}, \ldots, X_{nn}$ for which each $X_{ni}$ is distributed as $\mathrm{N}(0,1)$, and the covariances $\rho_{ij} = \mathrm{cov}(X_{ni}, X_{nj})$ satisfy $|\rho_{ij}| \le a_n$ for all $i, j$ such that $1 \le i, j \le n$ and $|i - j| \ge n^{\beta}/3$. Define $Y_{ni}$ by adding the sparse signal $I_{ni}\nu_1$

to $X_{ni}$, as at (3.12). Here, to delineate the way in which performance of the max classifier does not vary with dependence, we take $\nu_1$ to equal the value that $\nu$, in (3.12), would take in the absence of correlation: $\nu_1 = \sqrt{2r \log n}$. Compare (3.11).

THEOREM 3.9. *If $\beta \in (\frac{1}{2}, 1)$ and $r > (1 - \sqrt{1-\beta})^2$, then the probability that the max classifier correctly detects that a signal is present when it is, and the probability that the classifier does not detect a signal when it is not present, both converge to 1 uniformly in time-series in $\mathcal{X}(\beta, a_n)$.*

**4. Numerical properties.** We report here a simulation study assessing the performance of $\mathrm{hc}_n$. The analysis involved selecting pairs $(r, \beta)$ above the detection boundary, and generating samples from either population. To simulate a stationary Gaussian process $X_{n1}, \ldots, X_{nn}$ with the covariance function specified in (2.1), we used the method suggested by Wood and Chan [29].

Recall that $\kappa = \alpha_0/\alpha$, $r' = (1-\kappa)r$ and $\beta' = (1-\kappa)\beta$, and that we obtain the process $Y_{n1}, \ldots, Y_{nn}$ by adding to $X_{n1}, \ldots, X_{nn}$ a total of $n^{1-\beta'}$ signals, each of strength $\sqrt{2r' \log n}$. For the results reported here we selected two values of $(\beta, r)$, and for each, studied the performance of $\mathrm{hc}_n$ for different values of $\alpha_0$. Specifically, we chose $(\beta, r) = (0.6, 0.35)$ and $(\beta, r) = (0.75, 0.5)$; each is a vertical distance 0.25 above the detection boundary at (3.13). For each $(\beta, r)$ we treated $n = 2^{16}$ and $n = 2^{20}$, $\alpha = 0.5$ and $\alpha_0 = 0.05, 0.1, 0.15, 0.20$, so that $\kappa = 0.1, 0.2, 0.3, 0.4$. We chose $n$ large since the signals are highly sparse.

Our project was implemented in the following sequence of steps: (1) Generate the stationary Gaussian process $\{X_{n1}, \ldots, X_{nn}\}$ having the covariance function at (2.1). (2) Defining $K = n^{1-\beta'}$, generate $K$ variates from the uniform distribution on the unit interval, ordering them as $u_1 < u_2 < \cdots < u_K$; and generate a sequence $I_{ni}$ of zeros and ones, taking the value 1 at $i = \langle n u_i \rangle$ for $i = 1, 2, \ldots, K$, and the value 0 elsewhere. (Here, $\langle x \rangle$ denotes the integer nearest to $x$.) (3) Put $Y_{ni} = X_{ni} + \sqrt{2r' \log n} I_{ni}$ for $1 \le i \le n$. (4) Construct $\mathrm{hc}_n$ from the data $X_{ni}$ and $Y_{ni}$, where in the definition of $\mathrm{hc}_n$, $\mathcal{S}_n = (-t_n, t_n)$ with $t_n = \bar{\Phi}^{-1}(n^{\kappa-1})$. (5) Repeat steps (1)–(4) 100 independent times. The matlab code can be found at www.stat.purdue.edu/ jinj/Research/software/HallandJin.

For $(\beta, r) = (0.6, 0.35)$, if we decide to reject the null hypothesis of "no signal" when $n^{-\kappa/2}\mathrm{hc}_n \ge 2.2$, then, for $n = 2^{16}$, the empirical probability among 100 runs of committing a type I error equals 0.00 in the respective cases $\alpha_0 = 0.05$, 0.10, 0.15 and 0.20, and the corresponding empirical probabilities of committing a type II error are 0.00, 0.00, 0.00 and 0.06. Increasing $n$ to $2^{20}$, the respective type I error probabilities become 0.03, 0.00, 0.00 and 0.01, and type II error probabilities decrease to 0.00 in each case.

TABLE 1
*Mean and standard deviation (SD) of $n^{-\kappa/2}\,\mathrm{hc}_n$*

| $n$ | | $2^{10}$ | $2^{12}$ | $2^{14}$ | $2^{16}$ | $2^{18}$ | $2^{20}$ |
|---|---|---|---|---|---|---|---|
| $n^{-\kappa/2}\mathrm{hc}_n$ | Mean | 1.0273 | 1.0504 | 0.9162 | 0.8851 | 0.8467 | 0.7823 |
| | SD | 0.6234 | 0.4290 | 0.3810 | 0.4086 | 0.3702 | 0.3805 |

For $(\beta, r) = (0.75, 0.5)$, if we reject the null hypothesis when $n^{-\kappa/2}\mathrm{hc}_n \geq 1.7$, then, for $n = 2^{16}$, the type I error probabilities are 0.25, 0.02, 0.02 and 0.03, and the type II error probabilities are 0.04, 0.10, 0.21 and 0.41; they decrease to 0.17, 0.03, 0.04 and 0.00, and to 0.00, 0.00, 0.00 and 0.02, respectively, when $n = 2^{20}$.

We also investigated how closely the values of $\mathrm{hc}_n$ accord to our asymptotic analysis, which states that under the null hypothesis, $\mathrm{hc}_n = O_p(n^{(\kappa/2)+\varepsilon})$ and $n^{(\kappa/2)-\varepsilon} = O_p(\mathrm{hc}_n)$ for each $\varepsilon > 0$. We therefore conducted a simulation study where we fixed $(\alpha_0, \alpha) = (0.1, 0.5)$, entailing $\kappa = 0.2$, and took $\log_2 n$ in the range 10(2)20. For each $n$ we generated a stationary Gaussian process with the covariance structure specified in (2.1), and calculated $\mathrm{hc}_n$. We then repeated the simulations 100 independent times. The results are tabulated in Table 1, which shows that as $n$ increases from $10^3$ to $10^6$, the mean and standard deviation of $n^{-\kappa/2}\mathrm{hc}_n$ alter by only 24% and 39%, respectively, over the thousand-fold range of values of $n$.

## 5. Technical arguments.

5.1. *Proof of Theorem* 3.1. The case where $t$ lies in a bounded interval is straightforward to treat, and the situation of large negative $t$ can be addressed analogously to that for large positive $t$. Therefore we confine attention to $t \geq 1$.

Let $X$ and $Y$ be jointly normally distributed with zero means, unit variances and correlation coefficient $\rho$. Result (3.1) can be derived using the following lemma.

LEMMA 5.1. *If $t \to \infty$ and $\rho \to 1$ in such a manner that $t^2(1-\rho) \to 0$, then*

$$1 - P(X > t | Y > t) \sim \frac{t(1-\rho)^{1/2}}{\pi^{1/2}}. \tag{5.1}$$

5.2. *Proof of Theorem* 3.2. For reasons that preface the proof of Theorem 3.1, it suffices to treat the case $t \geq 1$.

LEMMA 5.2. *Let* $Z_{k1}, \dots, Z_{km_k}$, *for* $1 \le k \le \ell$, *denote normal random variables with zero means, unit variances and* $\mathrm{cov}(Z_{k1}, Z_{ki}) = \rho_{ki} \ge 0$ *for* $1 \le i \le m_k$. *Put* $m = m_1 + \cdots + m_\ell$ *and*

$$\rho_{\min} = \min\{\rho_{ki} : 1 \le i \le m_k,\ 1 \le k \le \ell\}.$$

*Then, uniformly in* $m \ge 1$ *and* $t \ge 1$,

$$\begin{aligned} &|1 - P\{I(Z_{ki} > t) = I(Z_{k1} > t) \text{ for } 1 \le i \le m_k \text{ and } 1 \le k \le \ell\}| \\ &\qquad \le Am(1 - \rho_{\min})^{1/2} e^{-t^2/2}, \end{aligned} \tag{5.2}$$

*where* $A$ *is a positive absolute constant.*

PROOF. Denote the left-hand side of (5.2) by LHS. Then,

$$\begin{aligned} 1 - \mathrm{LHS} &\le \sum_{k=1}^{\ell} \sum_{i=2}^{m_k} P\{I(Z_{ki} > t) \ne I(Z_{k1} > t)\} \\ &= 2 \sum_{k=1}^{\ell} \sum_{i=2}^{m_k} P(Z_{ki} > t, Z_{k1} \le t) \\ &= 2\{1 - \Phi(t)\} \sum_{k=1}^{\ell} \sum_{i=2}^{m_k} \{1 - P(Z_{ki} > t | Z_{k1} > t)\}. \end{aligned}$$

Result (5.2), and hence Lemma 5.2, follows from this bound and (5.1).

Next we derive Theorem 3.2. Given $\varepsilon \in (0, \lambda)$, partition the set of integers in the interval $[1, n^\lambda]$ into $n^{\lambda-\varepsilon}$ nonoverlapping subintervals, each of length $n^\varepsilon$, where the subintervals of integers are ordered as $\mathcal{I}_1, \dots, \mathcal{I}_{n^{\lambda-\varepsilon}}$ from left to right along the real line. (We shall omit integer-part notation, and also omit the straightforward treatment of the case where the last interval is a fragment, shorter than the other intervals.) Let $i_j$ denote the integer furthest to the left in $\mathcal{I}_j$. Then, using Lemma 5.2, we deduce that uniformly in $t \ge 1$,

$$\begin{aligned} &P\{I(X_{ni} > t) = I(X_{ni_j} > t) \text{ for } i \in \mathcal{I}_j \text{ and } 1 \le j \le n^{\lambda-\varepsilon}\} \\ &\qquad = 1 - O\{n^\lambda (n^{\varepsilon\alpha}/n^{\alpha_0})^{1/2} e^{-t^2/2}\} = 1 - O(n^{\lambda + (\varepsilon\alpha/2) - (\alpha_0/2)} e^{-t^2/2}). \end{aligned} \tag{5.3}$$

Suppose $0 < \lambda_1 < \lambda$, and that we have partitioned the interval $[1, n^\lambda]$ into $n^{\lambda-\lambda_1}$ nonoverlapping subintervals, each of length $n^{\lambda_1}$, arranged in order as $\mathcal{J}_1, \dots, \mathcal{J}_{n^{\lambda-\lambda_1}}$ from left to right; and that, with $i_j$ denoting the integer furthest to the left in $\mathcal{J}_j$,

$$\begin{aligned} &P\{I(X_{ni} > t) = I(X_{ni_j} > t) \text{ for } i \in \mathcal{J}_j \text{ and } 1 \le j \le n^{\lambda-\lambda_1}\} \\ &\qquad = 1 - O(n^\gamma e^{-t^2/2}), \end{aligned} \tag{5.4}$$

uniformly in $t \ge 1$, for a constant $\gamma > 0$. Result (5.3) is an instance of (5.4), with $\lambda_1 = \varepsilon$ and $\gamma = \lambda + (\varepsilon\alpha/2) - (\alpha_0/2)$. We shall argue by induction, from (5.3) via (5.4).

Consider the sparsely sampled time-series $X_{ni_1}, \dots, X_{ni_{n^{\lambda-\lambda_1}}}$. Partition the equally spaced integers $i_1, \dots, i_{n^{\lambda-\lambda_1}}$, defined immediately above (5.4), into consecutive blocks $\mathcal{K}_1, \dots, \mathcal{K}_{n^{\lambda-\lambda_1-\lambda_2}}$, each containing $n^{\lambda_2}$ integers, where $0 < \lambda_2 < \lambda - \lambda_1$. Let $k_\ell$ denote the integer furthest to the left in $\mathcal{K}_\ell$. Then, by Lemma 5.2, we have uniformly in $t \ge 1$,

$$\begin{aligned} &P\{I(X_{ni_j} > t) = I(X_{nk_\ell} > t) \text{ for } j \in \mathcal{K}_\ell \text{ and } 1 \le \ell \le n^{\lambda-\lambda_1-\lambda_2}\} \\ &\qquad = 1 - O\{n^{\lambda-\lambda_1}(n^{\lambda_2\alpha}/n^{\alpha_0})^{1/2} e^{-t^2/2}\} \\ &\qquad = 1 - O(n^{\lambda-\lambda_1+(\lambda_2\alpha/2)-(\alpha_0/2)} e^{-t^2/2}). \end{aligned} \tag{5.5}$$

For each $\ell$, put $\mathcal{L}_\ell = \bigcup_{j\in\mathcal{K}_\ell} \mathcal{J}_j$. This gives a partition of $[0, n^\lambda]$ into $n^{\lambda-\lambda_1-\lambda_2}$ disjoint subintervals $\mathcal{L}_\ell$, with $1 \le \ell \le n^{\lambda-\lambda_1-\lambda_2}$, each $\mathcal{L}_\ell$ containing $n^{\lambda_1+\lambda_2}$ consecutive integers, and for which, in view of (5.4) and (5.5), the version of (5.4) holds with $i_j$ replaced by $m_j$ (here denoting the least integer in $\mathcal{L}_j$), $\mathcal{J}_j$ replaced by $\mathcal{L}_j$, $n^{\lambda-\lambda_1}$ replaced by $n^{\lambda-\lambda_1-\lambda_2}$ and $\gamma$ replaced by

$$\gamma' = \max\{\gamma, \lambda - \lambda_1 + \tfrac{1}{2}(\lambda_2\alpha - \alpha_0)\}. \tag{5.6}$$

That is, uniformly in $t \ge 1$,

$$\begin{aligned} &P\{I(X_{ni} > t) = I(X_{nm_j} > t) \text{ for } i \in \mathcal{L}_j \text{ and } 1 \le j \le n^{\lambda-\lambda_1-\lambda_2}\} \\ &\qquad = 1 - O(n^{\gamma'} e^{-t^2/2}). \end{aligned} \tag{5.7}$$

Having achieved, in (5.3), result (5.4) for $\lambda_1 = \varepsilon$ and $\gamma = \lambda + \frac{1}{2}(\varepsilon\alpha - \alpha_0)$, we may, noting the definition of $\gamma'$ at (5.6), ensure that (5.4) holds with $\gamma' = \gamma$, by choosing $\lambda_2$ such that

$$\lambda - \lambda_1 + \tfrac{1}{2}(\lambda_2\alpha - \alpha_0) = \lambda + \tfrac{1}{2}(\varepsilon\alpha - \alpha_0);$$

that is, $\lambda_2 = \varepsilon + (2\lambda_1/\alpha)$. Selecting this $\lambda_2$ implies that, in passing from (5.3) to (5.7), the number of subintervals [$\mathcal{I}_j$ in (5.3) and $\mathcal{L}_j$ in (5.7)] has decreased from $n^{\lambda-\varepsilon}$ to $n^{\xi_1}$, where

$$\xi_1 = \lambda - \lambda_1 - \lambda_2 = \lambda - \varepsilon - \lambda_1\left(1 + \frac{2}{\alpha}\right) = \lambda - \lambda_1',$$

with $\lambda_1' = \varepsilon + \lambda_1(1 + 2\alpha^{-1}) = \varepsilon(1+\tau)$ and $\tau = 1 + 2\alpha^{-1}$. That is, (5.4) holds with $\lambda_1$ replaced by $\lambda_1'$ but still with $\gamma = \lambda + \frac{1}{2}(\varepsilon\alpha - \alpha_0)$. A further iteration of this argument decreases $\xi_1$ to

$$\xi_2 = \lambda - \varepsilon - \varepsilon\left(1 + \frac{2}{\alpha}\right) - \left[\varepsilon + \frac{2}{\alpha}\left\{\varepsilon + \varepsilon\left(1 + \frac{2}{\alpha}\right)\right\}\right] = \lambda - \varepsilon(1 + \tau + \tau^2),$$

achieving $\xi_k = \lambda - \varepsilon(1+\tau+\cdots+\tau^k)$ after a further $k-2$ iterations. Stopping when $\xi_k < 0$, we conclude that, uniformly in $t \geq 1$,

$$
\begin{aligned}
&P\{I(X_{ni} > t) = I(X_{n1} > t) \text{ for } 1 \leq i \leq n^{\lambda}\} \\
&\quad = 1 - O(n^{\lambda+(\varepsilon\alpha/2)-(\alpha_0/2)} e^{-t^2/2}).
\end{aligned} \tag{5.8}
$$

Since $\varepsilon > 0$ is arbitrary, then (5.8) implies (3.3), establishing Theorem 3.2. □

5.3. *Proof of Theorem* 3.3. It suffices to treat the case $t \geq 1$. Theorem 3.2 implies that, uniformly in $t \geq 1$,

$$
\max_{1\leq j\leq b} P\{I(X_{ni} > t) = I(X_{ns_j} > t) \text{ for } i \in \mathcal{B}_j\} = 1 - O(n^{\lambda-(\alpha_0/2)+\eta} e^{-t^2/2}).
$$

Therefore, the probability on the left-hand side of (3.7) equals

$$
1 - O(bn^{\lambda-(\alpha_0/2)+\eta} e^{-t^2/2}),
$$

uniformly in $t$. Theorem 3.3 follows from this property.

5.4. *Proof of Theorem* 3.4. Let $Z$ have the $\mathrm{N}(0,1)$ distribution, and choose $\lambda < \kappa$ so close to $\kappa$ that $1 - \lambda - q < 0$. Write LHS for the left-hand side of (3.8). Then, using Theorem 3.3 to derive the inequality below, we have

$$
\begin{aligned}
\mathrm{LHS} &\leq \sum_{j=1}^{b} P(X_{ns_j} > t) + o(1) = bP(Z > t) + o(1) \\
&= O(n^{1-\lambda} e^{-t^2/2}) + o(1) = O(n^{1-\lambda-q}) + o(1) \to 0,
\end{aligned}
$$

which implies (3.8).

5.5. *Proof of Theorem* 3.5. Theorem 3.4 implies that, unless $q \leq 1 - \kappa$, the level $t = \sqrt{2q \log n}$ is too high if we are conducting inference for the process $X_{n1}, \ldots, X_{nn}$. Therefore we take $q < 1 - \kappa$ in Lemma 5.3, below.

LEMMA 5.3. *Assume that the stationary, zero-mean Gaussian process $X_{ni}$ has autocovariance given by* (2.1), *and that $\kappa < 1$. Let $t_n$ denote a sequence of positive constants diverging to infinity in such a manner that $n^{1-\kappa}\bar{\Phi}(t_n)$ is bounded away from zero. Then, for each integer $\nu \geq 1$,*

$$
E\left[\sum_{i=1}^{n} \{I(X_{ni} > t) - P(X_{ni} > t)\}\right]^{2\nu} = O[\{n^{\kappa+1} \Phi(t) \bar{\Phi}(t)\}^{\nu}], \tag{5.9}
$$

*uniformly in $|t| \leq t_n$.*

To complete the proof of Theorem 3.5, define

$$\delta(t) = \{n^{\kappa+1}\Phi(t)\bar{\Phi}(t)\}^{1/2},$$

$$\Delta(t) = \sum_{i=1}^{n}\{I(X_{ni} > t) - P(X_{ni} > t)\},$$

and observe that, if $\mathcal{T}_n$ is any subset of the values of $t$ in the interval $[-t_n, t_n]$, and if $\eta > 0$, then by Markov's inequality,

$$\begin{aligned} P\{|\Delta(t)| > n^{\eta}\delta(t) \text{ for some } t \in \mathcal{T}_n\} &\le n^{-2\nu\eta} E\Big\{\sup_{t\in\mathcal{T}_n} |\Delta(t)/\delta(t)|^{2\nu}\Big\} \\ &\le n^{-2\nu\eta}\sum_{t\in\mathcal{T}_n} E\{|\Delta(t)/\delta(t)|^{2\nu}\} \\ &\le n^{-2\nu\eta}(\#\mathcal{T}_n)\sup_{t\in\mathcal{T}_n} E\{|\Delta(t)/\delta(t)|^{2\nu}\} \\ &= O\{n^{-2\nu\eta}(\#\mathcal{T}_n)\}, \end{aligned} \tag{5.10}$$

where the last identity follows from Lemma 5.3. Therefore, as long as $\mathcal{T}_n$ has no more than polynomially many elements, the following result holds: For all constants $C_1, \eta > 0$,

$$P\{|\Delta(t)| > n^{\eta}\delta(t) \text{ for some } t \in \mathcal{T}_n\} = O(n^{-C_1}). \tag{5.11}$$

By choosing the elements of $\mathcal{T}_n$ to be equally spaced on $[-t_n, t_n]$, and selecting the spacing to equal $n^{-C_2}$, for $C_2$ sufficiently large but fixed; and noting that $\kappa + 1 - q > 0$ and, uniformly in $t \in [-t_n, t_n]$,

$$\delta(t) \ge \{n^{\kappa+1-q}(\log n)^{-1}\}^{1/2};$$

we deduce from (5.11) that

$$P\{|\Delta(t)| > n^{\eta}\delta(t) \text{ for some } t \in [-t_n, t_n]\} \to 0. \tag{5.12}$$

This implies (3.9).

Proof of Theorem 3.7. Recall that we distribute $n^{\xi} = N^{1-\beta}$ signals among the $n$ time-series points, where $\xi = (1-\beta)(1-\kappa)$. Let

$$U_n(t) = \sum_{i=1}^{n}\{I(X_{ni} > t) - P(X_{ni} > t)\},$$

$$V_n(t) = \sum_{i=1}^{n}\{I(Y_{ni} > t) - P(Y_{ni} > t)\},$$

$$W_n(t) = \sum_{i=1}^{n} I_{ni}[I(Y_{ni} > t) - I(X_{ni} > t) - \{P(Y_{ni} > t) - P(X_{ni} > t)\}],$$

where $I_{ni} = 1$ if a signal is added at "time" point $i$, equaling zero if no signal is added there. Define $\Psi_1(t) = \Phi(t) - \Phi(t-\nu)$ and $\Psi_2 = \Psi_1(1-\Psi_1)$. An argument similar to that used to derive Lemma 5.3 may be employed to show that, and uniformly in $|t| \le t_n$,

$$E\{W_n(t)^{2\nu}\} = O[n^{2\nu\kappa}\{n^{\xi}\Psi_2(t)\}^{\nu}]. \tag{5.13}$$

The right-hand side of (5.13) is an upper bound, over all choices of the distribution of signal among the time-series data $X_{ni}$ as well as over all values of $|t| \le t_n$. The order of magnitude of the left-hand side is maximized when just $N^{1-\beta}$ of the consecutive blocks of indices $\{1,\dots,n^{\kappa}\}, \{n^{\kappa}+1,\dots,2n^{\kappa}\},\dots$ are chosen to receive signals, and, for each of the chosen blocks, a signal is applied to each value of $X_{ni}$ the index of which is in that block.

Result (5.13), Markov's inequality and the argument leading to (5.10) imply that if $\mathcal{T}_n$ is any finite subset of values in the interval $[-t_n, t_n]$, then for each $\eta > 0$,

$$\begin{aligned} &P[|W_n(t)| > n^{\kappa+\eta}\{n^{\xi}\Psi_2(t)\}^{1/2} \text{ for some } t \in \mathcal{T}_n] \\ &\qquad \le n^{-2\nu(\kappa+\eta)} E\left[\sup_{t\in\mathcal{T}_n}\left|\frac{W_n(t)}{\{n^{\xi}\Psi_2(t)\}^{1/2}}\right|^{2\nu}\right] = O\{n^{-2\nu\eta}(\#\mathcal{T}_n)\}. \end{aligned}$$

This leads to the following analogue of (5.12): for each $\eta > 0$,

$$P[|W_n(t)| > n^{\kappa+\eta}\{n^{\xi}\Psi_2(t)\}^{1/2} \text{ for some } t \in [-t_n, t_n]] \to 0. \tag{5.14}$$

Noting that $V_n = U_n + W_n$, writing $\Psi_3 = \Phi\bar{\Phi}$, and combining (5.11) and (5.14), we deduce that for each $\eta > 0$,

$$\begin{aligned} &P\left[\frac{|V_n(t)|}{\{N\Psi_3(t)\}^{1/2}} > n^{\eta}\left\{n^{\xi+3\kappa-1}\frac{\Psi_2(t)}{\Psi_3(t)} + n^{2\kappa}\right\}^{1/2}\right. \\ &\qquad\qquad\qquad \left.\text{for some } t \in [-t_n, t_n]\vphantom{\frac{1}{1}}\right] \to 0. \end{aligned} \tag{5.15}$$

Note too that

$$\begin{aligned} u_n(t) &\equiv \frac{\sum_i I_{ni} E\{I(Y_{ni} > t) - I(X_{ni} > t)\}}{\{N\Phi_3(t)\}^{1/2}} \\ &= N^{1-\beta} n^{\kappa} \frac{\Phi(t) - \Phi(t-\nu)}{\{N\Phi_3(t)\}^{1/2}} = \left\{n^{2\xi+3\kappa-1}\frac{\Psi_1(t)^2}{\Psi_3(t)}\right\}^{1/2}. \end{aligned}$$

Result (3.14) will follow from this property and (5.15), provided we show that, for some $\eta > 0$ and uniformly in $t \in [-t_n, t_n]$,

$$n^{\xi+3\kappa-1}\frac{\Psi_1(t)}{\Psi_3(t)} + n^{2\kappa} = O\left\{n^{2\xi+3\kappa-1-\eta}\frac{\Psi_1(t)^2}{\Psi_3(t)}\right\},$$

or equivalently, for some $\eta > 0$ and uniformly in $t \in [-t_n, t_n]$,

$$n^{\eta/2} = O\{n^{\xi}\Psi_1(t)\}, \qquad n^{\eta} = O\{n^{2\xi+\kappa-1}\Psi_1(t)^2\Psi_3(t)^{-1}\}. \tag{5.16}$$

Observe that $n^{1-\kappa}\Psi_3(t)$ is bounded away from zero uniformly in $|t| \le t_n$, and so the first part of (5.16) holds provided the second part does. Now, the second part is equivalent to $n^{\eta} = O[\{n^{-\kappa}u_n(t)\}^2]$ uniformly in $t$. If the point $(\beta, r)$ lies an amount $\varepsilon > 0$ above the standard detection boundary, at (3.13), then there exists $\delta = \delta(\beta, \varepsilon) > 0$, and $t = t(n) \in [-t_n, t_n]$, such that $n^{-\kappa}u_n(t) \ge \text{const.}\ n^{\delta}$ for all $n$. Therefore, the second part of (5.16) holds if $0 < \eta < \delta(\beta, \varepsilon)$. □

PROOF OF THEOREM 3.8. Note that, for each $i$, $X_{n,i+1} - X_{ni}$ is normally distributed with zero mean and variance $2/n^{\alpha_0}$. Therefore, the probability that the signal detector determines that the signal is present, given that it is not, is bounded above by

$$2n[1 - \Phi\{(2C\log n)^{1/2}\}] = O\{n^{1-C}(\log n)^{-1/2}\} \to 0;$$

and the probability that the signal is found by the detector to be present, given that it is, is bounded below by

$$1 - \Phi\{(2C\log n)^{1/2} - (2r\log n)^{1/2}\} \to 1. \qquad \square$$

Department of Mathematics and Statistics
University of Melbourne
Melbourne, VIC 3010
Australia
and
Department of Statistics
University of California, Davis
One Shields Avenue
Davis, California 95616
USA

Department of Statistics
Purdue University
250 North University Street
West Lafayette, Indiana 47907–2066
USA
E-mail: jinj@stat.purdue.edu